\documentclass[a4paper,10pt]{article}

\def \Z {{\mathbf {Z}}}

\def \N {{\mathbf {N}}}

\def\u{\bigsqcup}
\def\eps{\varepsilon}
\title{Weak Closure Theorem for Double Staircase Actions}
\author{V.V. Ryzhikov\footnote{\large This work is partially  supported by  the grant
NSh 8508.2010.1.}}

\begin{document}
\Large
\maketitle

\section{Introduction} 
Considering  measure-preserving transformations   of a Probability space $(X,\mu$)
we introduce  a double staircase  construction $T$ and show that its
semi-group of all weak limits of powers ($WLP(T)$) is 
$$\{\Theta, 2^{-m}T^n+(1-2^{-m})\Theta \ :  m\in\N,\  n\in \Z\},$$
where $\Theta$ stands for  the orthogonal projection $L_2(X,\mu)$ onto the space of constant functions.
Mixing sequences are controlled  via Adams'  method \cite{A}, non-mixing ones  come  to light  
 by means of   secondary  limit methods (see \cite{R}). 
\\
{\bf Staircase rank one transformation} is determined by an integer $h_1$ and a sequence $r_j$  of cuttings.
  We recall its definition.
Let our  $T$ on the step $j$  be associated with  a collection of disjoint intervals
$$E_j, TE_j, T^2E_j,\dots, T^{h_j}E_j.$$
We cut $E_j$ into $r_j$ subintervals  of the same measure
$$E_j=E_j^1\u E_j^2\u  E_j^3\u\dots\u E_j^{r_j},$$  
then for all $i=1,2,\dots, r_j$ we  consider columns
$$E_j^i, TE_j^i ,T^2 E_j^i,\dots, T^{h_j}E_j^i.$$
Adding  over $i$-th column $i-1$ spacers   we obtain a partition  
\vspace{4mm}

$E_j^1, TE_j^1, T^2 E_j^1,\dots, T^{h_j-1}E_j^1,T^{h_j}E_j^1,$

$E_j^2, TE_j^2 ,T^2 E_j^2,\dots, T^{h_j-1}E_j^2,T^{h_j}E_j^2,T^{h_j+1}E_j^2,$

$E_j^3, TE_j^3 ,T^2 E_j^3,\dots,T^{h_j-1}E_j^3, T^{h_j}E_j^3,T^{h_j+1}E_j^3,T^{h_j+2}E_j^3,$

$\dots\ \ \dots \ \ \dots$

$E_j^r, TE_j^r ,T^2 E_j^r,\dots,T^{h_j-1}E_j^r, T^{h_j}E_j^r,T^{h_j+1}E_j^r,T^{h_j+2}E_j^r,\dots, T^{h_j+r-1}E_j^r,$
\vspace{2mm}

where   $r=r_j$.    For all  $i<r_j$ we set
$T^{h_j+i}E_j^i = E_j^{i+1}.$
\\ Thus,  we get  $j+1$-the tower 
$E_{j+1}, TE_{j+1} T^2 E_{j+1},\dots, T^{h_{j+1}}E_{j+1}$
with
 $$E_{j+1}= E^1_j, \ \
h_{j+1}+1=(h_j+1)r_j +\sum_{i=1}^{r_j-1}i.$$

This staircase construction is a special case ($s(i)=i-1$)  of a general rank one construction
with  a sequence $\bar s_j$ of spacer vectors
$$ \bar s_j=(s_j(1), s_j(2),\dots, s_j(r_j-1),s_j(r_j)).$$  
\\
 {\bf    Double staircase  transformation} is defined by an integer $h_1$ and a   spacer sequence 
$$ \bar s_j=(0,1,2,3,\dots, r'_j-2,r'_j-1, 0,1,2,3,\dots, r'_j-2,r'_j-1),$$
where $2r'_j=r_j$. We  presume that $r_j\to\infty$.
In what follows we presume that  Adams' condition $r_j^2/h_j\to0$ is  satisfied.   This restriction plays only a  technical role 
and will be used implicitly  in approximation formulas.
\\
{\bf On notations.}  We write $a(j)\approx b(j)$ instead of $a(j)-b(j)\to 0$ (or $\frac{a(j)}{b(j)}\to 1$) and  use weak $\approx_w$ and strong $\approx$ operator approximations.
\\
{\bf Main result.}

{\bf THEOREM 1.} {\it   A double staircase transformation $T$  possesses the following semi-group  of weak limits of 
its powers ($WLP$):  $$\{\Theta, 2^{-m}T^n+(1-2^{-m})\Theta \ :  m\in\N,\  n\in \Z\},$$
where $\Theta$ is the orthogonal projection into the space of constant functions.}

{}

\bf LEMMA 1. \it  $M=\N\setminus \bigcup_j [0.5h_{j+1}-h_j, 0.5h_{j+1}+h_j]$ is a mixing set:
if $m_i\to\infty$ and $m_i\in M$, then $T^{m_i}\to\Theta$.
\vspace{3mm}
\rm 

Thus,  if   $n_j$ is a non-mixing sequence, $|n_j|\in[h_{j'},{h_{j'+1}})$, then $$|n_j|\in  [0.5h_{j'+1}-h_{j'}, 0.5h_{j'+1}+h_{j'}].$$

We prove Lemma 1 in section 3.
\section{ Non-mixing sequences} 
 An operator of multiplication by $\chi_D$ is denoted by  $\hat D$. 
\vspace{2mm}\\
 \bf LEMMA 2. \it  Let $\mu(D(j))\to a>0$ and $\mu(D(j)\Delta TD(j))\to 0$.  Then  $\hat D(j) \approx_w aI $.
If $Q(j)\approx\Theta$, then 
$\hat D(j)Q(j)\approx_w a\Theta$.  \rm
\vspace{2mm} 

A proof of this lemma is an exercise.     
\vspace{2mm}
\\
\bf Example 1.  \rm Let, for instance, $n_j=h_{j+1}/2,$  then  standard rank one calculations  (see \cite{A},\cite{R}) allow us to

$\mu(T^{n_j}A\cap B)\approx $
 $$\approx \frac{1}{r_{j}}\sum_{i=1}^{r_{j}/2}\mu(T^{-d_0i}A\cap B)+ 
\frac{1}{r_{j+1}}\sum_{i=1}^{r_{j+1}}\mu(T^{-h_{j+1}/2}T^{-d_1i}A\cap B\cap  D_0(j)), \eqno (\ast)
$$
where $d_0=0$, $d_1=1$,  $$D_0(j)=\u_{i=0}^{h_{j+1}/2-1} E_{j+1}, \ \
\mu(D_0)\approx \frac{1}{2}.$$
 We rewrite $(\ast)$ in the form
$$T^{n_j}\approx_w \frac{1}{2}I+\hat D_0(j)\frac{1}{r_{j+1}}\sum_{i=1}^{r_{j+1}}T^{-h_{j+1}/2}T^{-i}.$$
From the ergodicity of $T$ we get 
  $$
\frac{1}{r_{j+1}}\sum_{i=1}^{r_{j+1}}T^{-h_{j+1}/2}T^{-i}\approx \Theta,$$
hence,
$$T^{n_j}\approx_w  \frac{1}{2}I+\frac{1}{2}\Theta.$$

\bf Example 2.  \rm Let   $n_j=h_{j+1}/2+k_j, \ \ 0\leq |k_j| \leq h_j,$ \\ we have
$$\mu(T^{n_j}A\cap B)\approx \mu(T^{k_j}A\cap B\cap D(j))+
\frac{1}{r_{j}}\sum_{i=1}^{r_{j}}\mu(T^{k(1,j)}T^{-i}A\cap B\cap D_1(j))+$$
$$+\frac{1}{r_{j+1}}\sum_{i=1}^{r_{j+1}}\mu(T^{-h_{j+1}+n_j}T^{-i}A\cap B\cap D_0(j)),
$$
where  $$D(j)=U_j \cap\u_{i=|k_j|}^{h_{j}}T^i E_{j}, \ \ D_1(j)=U_j\cap \u^{|k_j|-1}_{i=0} T^iE_{j}, \ \ U_j=\u_{i=h_{j+1}/2}^{h_{j+1}}T^i E_{j+1}.$$
From the ergodicity of $T$  we get
$$\mu(T^{n_j}A\cap B)\approx \mu(T^{k_j}A\cap B\cap D(j)) +\left(\frac{h_j-|k_j|}{2h_j}+\frac{1}{2}\right)
\mu(A)\mu( B),$$
$$T^{n_j}\approx_w \hat D(j)T^{k_j} +\left(\frac{h_j-|k_j|}{2h_j}+\frac{1}{2}\right)\Theta.$$
\\
 If $n_j$ is not mixing, \rm then  from Lemma 2 it follows that $k_j$ is not mixing too. We get
  the following alternative \it

   either
$$T^{n_j}\approx_w \hat D(j)T^{k_j} +\frac{3}{4}\Theta, \ \ \mu(D(j))\approx  \frac{1}{4},\  \ |k_j|/h_j \approx 1/2,$$

or
$$T^{n_j}\approx_w \hat D(j)T^{k_j} +\frac{1}{2}\Theta, \ \ \mu(D(j))\approx  \frac{1}{2}, \ \ k_j/h_j \approx 0.$$ \rm
\\
\bf In fact \rm we have a strong approximation
$$T^{n_j}\approx \hat D(j)T^{k_j} +\hat Y(j)T^{n_j},$$
where  $Y(j)=X\setminus D(j)$ and   $$\hat Y(j)T^{n_j}\approx_w (1-\mu(D(j))\Theta.$$
\bf Let $k_j$ be not bounded, \rm then  there is $k'_j$,  $k'_j<< k_j$  such that
$$T^{k_j}\approx \hat D'(j)T^{k'_j} +\hat Y'(j)T^{n_j},$$
$$T^{n_j}\approx \hat D(j)(\hat D'(j)T^{k'_j} +\hat Y'(j)T^{k_j}) +\hat Y(j)T^{n_j}.$$
\bf If $k'_j$ is not bounded,  \rm then for $k_j''<<k_j'$ we have
$$T^{n_j}\approx \hat D(j)\left[\hat D'(j)\left(\hat D''(j)T^{k''_j} +\hat Y''(j)T^{k'_j} \right)+\hat Y'(j)T^{k_j}\right]+\hat Y(j)T^{n_j}.$$

Omitting $(j)$ in $D(j),\dots Y''(j)$ we rewrite the latter by 
$$T^{n_j}\approx \hat D\hat D'\hat D''T^{k''_j} + \hat D\hat D'Y''T^{k'_j}+
\hat D\hat Y'T^{k_j}+ \hat YT^{n_j}.$$
And so on.
However the number of  iterations $'$  must be finite. Indeed,  $D, D', D''$ are special pieces  of different towers, so
$$\mu(D\cap D'\cap D'')\approx \mu(D)\mu(D')\mu(D'')\approx \frac{1}{2^m}, \ 3
\leq m\leq 6$$
 (they are   almost independent), and 
$$\hat YT^{n_j}\approx_w \mu(Y)\Theta,$$
$$\hat D\hat Y'T^{k_j}\approx_w \mu(D)\mu(Y')\Theta$$
(this follows from $\hat Y'T^{k_j}\approx_w \mu(Y')\Theta$ and $\mu (D\cap Y')\approx \mu (D)\mu (Y')$),
$$\hat D\hat D'Y''T^{k'_j}\approx_w \mu(D)\mu(D')\mu(Y'')\Theta.$$

If  all   sequences $k^{''\dots '}_j$ are unbounded, then 
$$ \mu(D\cap D'\dots \cap D''^\dots)\approx 0,$$
 $$T^{n_j}\to \Theta,$$  a contradiction ($n_j$ is non-mixing).  Thus, there is $k^{''\dots '}_j=k$, and  $m$ such that  
$$\mu\left(D\cap D'\cap\dots \cap D^{''\dots'}\right)\to  \frac{1}{2^{m}}, \ \ so,\ \ \hat D\hat  D'\dots \hat D^{''\dots'}T^k\to_w \frac{1}{2^{m}}T^k,$$
$$T^{n_j}\to \frac{1}{2^{m}}T^k+ \left(1-\frac{1}{2^{m}}\right)\Theta.$$

\section{ Mixing sequences.}
Now we  prove Lemma 1.  To a reader who is familiar with Adams' approach we can explain a proof  in "two
words".   Adams  found a total control of mixing for staircase transformation based
on mixing sequences  $m_j\in [h_j, Ch_j]$.  For this  he produced a non-trivial method to control  mixing properties ($P(j)\approx\Theta$) of averaging operators in a form
$$P(j)=\frac{1}{r(j)}\sum_{i=1}^{r(j)}  T^{d(j)i},$$  
where some special sequences $r(j), d(j)$ ($r(j)\to\infty$).      The case $C>d(m)>0$ is trivial  ($T$ is totally ergodic);
the case $d(j) \to \infty$ is of interest.  
Adams found a special number $q(m)$ such that $q(j)<<r(j)$ and for large $L$ ($L=L(j)$ tends to infinity very slowly)
$$A(j)=\frac{1}{L}\sum_{i=1}^{L}  T^{q(j)d(j)i}\approx \Theta.$$  
Then for all $n$ 
$$T^nA(j)\approx \Theta$$
holds, so
$$P(j)\approx\Theta$$
as a convex sum of  $T^nA(j)$.   
The principle personage $q(j)$ has been  extremely resourcefully selected by Adams: he found $q(j)$,
$h_{p(j)}< q(j)d(j) <2h_{p(j)}$, such that  
$$T^{h(j)d(j)},T^{2h(j)d(j)},\dots,T^{Lh(j)d(j)}\approx_w \Theta.$$
Hence,  
$$A(j)^\ast A(j)\approx_w\Theta,$$
$$A(j)\approx\Theta,  \ \ P(j)\approx\Theta. $$

What is changed in our double staircase situations?  Almost nothing.  Again we have a starting mixing
sequence $[h_j, Ch_j]$  ($1<< C $). Again we can approximate $T^{m_j}$ by similar averaging operators
(now we deal with 5 approximating operators  instead of Adams' 3 operators, but all operators  are of the same nature).
All is similar except  one  thing: now the case
$$d(m)=0$$ 
 may appear.  Here we have
the following effect:  the image $T^mE_j$ of our the base $E_j$ has a flat part that  is situated in
one of floors $T^kE_j$, $0\leq k\leq h_j$. 
For example,  the case  $n_j=h_{j+1}/2$ (it has been considered above) gives
$$2\mu(T^{n_j}A\cap B)\approx \mu(A\cap B) +\frac{1}{r_{j+1}}\sum_{i=1}^{r_{j+1}}\mu(T^{-i}A\cap B)= $$  
 $$ = \frac{2}{r_{j}}\sum_{i=1}^{r_{j}/2}\mu(T^{-di}A\cap B)+ \frac{1}{r_{j+1}}\sum_{i=1}^{r_{j+1}}\mu(T^{-h_{j+1}/2}T^{-d_0i}A\cap B),$$
where $d=0$, $d_0=1$. (A half of the image $T^mE_j$   is  in
 $E_j$.)

For $n_j=h_{j+1}+h_j$ we now have $d=1$, $d_0=1$,
$$T^{n_j}\approx_w \frac{1}{r_{j}}\sum_{i=1}^{r_{j}/2}T^{-di}+\frac{1}{2r_{j+1}}\sum_{i=1}^{r_{j+1}}T^{-h_{j+1}/2-h_j}T^{-d_0i}\approx_w \Theta.\eqno (1)$$

\bf Generally \rm we get an approximation 
$$T^{m_j}\approx_w \sum_{s=0}^4   \hat D_s(j) T^{k(s,j)}P_s(j),$$ 
where  $D_s(j)$  tile together the space $X$,  $|k(s,j)|\leq h_j$,
$$P_s(j)=\frac{1}{r(s,j)}\sum_{i=1}^{r(s,j)}  T^{d(s,j)i}.$$
Let for instant  $m_j /h_{j+1}\approx 0.2$, then we have the following situation:
$$r(0,j)=h_{j+1}-1, \  r(1,j)=r(2,j)\approx 0.3r_j, \  r(3,j)=r(4,j)\approx 0.2r_j;$$
$$\mu(D_0(j))\approx 0.2,\ \mu(D_1(j))+ \mu(D_2(j)) \approx 0.6, \ \mu(D_3(j)) + \mu(D_4(j))
\approx 0.2;$$
$$d(0,j)=1,  \ d(1,j)=d(2,j)+1 \approx 0,3r_j, \ \ d(3,j)+1=d(4,j) \approx - 0,1r_j.$$

Sometimes certain $D_s(j)$ could be  vanishing.  For example, in our formula (1) (as $m_j=0.5h_{j+1}+h_j$) we had $\mu(D_0(j))=0.5$, $\mu(D_4(j))=0.5$,  $d(0,j)=d_0=1=d=d(4,j)$ (here $D_0(j)$ is a part of $j$-tower that
is situated under the first stairs array of spacers,  $D_4(j)$ -- under the second one). 
We will not weary the reader with tedious calculations,  we may  claim: for 
$$m_j \in [h_j, h_{j+1}]\setminus 
[h_{j+1}/2-h_j, h_{j+1}/2+h_j]$$ the corresponding    $d(s,j)$ are non-zero.  Thus,  Adams' method guarantees our $\{m_j\}$ to be mixing:
$$T^{m_j}\approx_w \sum_{s=0}^4   \hat D_s(j) T^{k(s,j)}P_s(j)\approx \Theta.$$

 Lemma 1 is proved. Theorem 1 is proved. 
 
\section{Remarks on related infinite  transformations}

1. There is a simple approach  to construct non-mixing Gaussian  automorphisms (see also \cite{DP}, p. 92.)  and  Poisson suspensions    
with explicit countable WLP.  Following \cite{R1} let's consider  double Sidon rank one infinite transformation $T$:  a rank one construction with a double spacer sequence 
$$s_j(1),s_j(2),\dots, s_j(r'_j),s_j(1),s_j(2),\dots, s_j(r'_j) $$  satisfied the condition
$$ h_j<<s_j(1)<< s_j(2)<<\dots<< s_j(r'_j-1)<<s_j(r'_j).$$ 
Easily controlling mixing sequences  we  get  simply $WLP(T)=\{0, 2^{-m}T^n\}$ (let us remark also that the centralizer of $T$ is  trivial), so  
$$WLP({\bf exp(}T{\bf )})=\{\Theta, {\bf exp (}2^{-m}T^n{\bf )}:  m=1, 2, \dots, \  n\in\Z\}.$$ \rm

The work \cite{DR} gives also a natural way to construct a double spacer map $S$ with  $WLP(S)=\{0, 2^{-m}S^n\}$.

2. 
A  calculation WLP(T) for  double staircase transformations $T$ of an infinite measure space  is  an interesting (and maybe hard) problem  as $r_j \sim h_j$.  

\it Conjecture.  Any staircase transformation   of an  infinite measure space is mixing as  $r_j\to\infty$.\rm

The  case $r_j/h_j\to 0$ is solved  ( the author presented a proof at Roscoff conference 
 "Stochastic properties of dynamical systems and random walks", June, 2010).

\bf  Problem. \it    Prove the mixing for the "$r_j=h_j$" staircase transformation.\rm
\\
The calculation    of $|P\cap P+m|$ for  
$$P=\{p\ : \ p=dh+\sum_{i=0}^{d-1} i, \ d=1,2,\dots, [(1-\eps)h]\}$$
 is naturally connected to the problem and  seems   non-trivial.

 E-mail: vryz@mail.ru
\end{document}